\documentclass[graybox]{svmult}

\usepackage{mathptmx}       
\usepackage{helvet}         
\usepackage{courier}        
\usepackage{type1cm}        
\usepackage{makeidx}         
\usepackage{graphicx}        
\usepackage{multicol}        
\usepackage[bottom]{footmisc}

\usepackage{amsmath}
\usepackage{amssymb}
\usepackage{amscd}
\usepackage{indentfirst}

\def\1{{\bf 1}}

\makeindex             

\usepackage[T1]{fontenc}                         
\usepackage{amsmath,amsfonts}                    
\DeclareMathAlphabet{\mathbbmsl}{U}{bbm}{m}{sl}  

\begin{document}

\title*{Double-sided Taylor's approximations and their applications in Theory of analytic inequalities}
\titlerunning{Double-sided Taylor's approximations and their app. in Theory of analytic inequalities}
\author{Branko Male\v sevi\' c, Marija Ra\v sajski, Tatjana Lutovac}
\authorrunning{B.~Male\v sevi\' c, M.~Ra\v sajski, T.~Lutovac}
\institute{B.~Male\v sevi\' c, M.~Ra\v sajski, T.~Lutovac \at School of Electrical Engineering, University of Belgrade\\
\email{branko.malesevic@etf.bg.ac.rs}\\
\email{marija.rasajski@etf.bg.ac.rs} \\
\email{tatjana.lutovac@etf.bg.ac.rs}}

\maketitle

\vspace*{-10.0 mm}

\abstract{In this paper the double-sided {\sc Taylor}'s approximations are studied. A short proof of a well-known theorem on the double-sided
{\sc Taylor}'s approximations is introduced. Also, two new theorems are proved regarding the monotonicity of such approximations. Then we present some new applications of the double-sided {\sc Taylor}'s approximations in the theory of analytic inequalities.}

\section{Introduction}

Consider a real function $f : (a, b) \longrightarrow \mathbb{R}$ such that there exist finite limits
 $f^{(k)}(a+)=\!\lim\limits_{x \rightarrow a+}{f^{(k)}(x)}$, for $k=0,1,\ldots,n$.
Let us denote by $T_n^{f,\,a+}(x)$ {\sc Taylor}'s polynomial of degree $n$, $n \!\in\! \mathbb{N}_{0}$, for the function $f(x)$ in the right neighborhood~of~$a$:
$$
T_n^{f,\,a+}(x)=\displaystyle\sum_{k=0}^{n}{\displaystyle\frac{f^{(k)}(a+)}{k!}(x-a)^k}.
$$
We will call $T_n^{f,\,a+}(x)$ the {\em first {\sc Taylor}'s approximation in the right neighborhood~of~$a$}.

\medskip
Similarly,  the {\em first {\sc Taylor}'s approximation in the left neighborhood of $b$} is defined by:
$$
T_n^{f,\,b-}(x)=\displaystyle\sum_{k=0}^{n}{\displaystyle\frac{f^{(k)}(b-)}{k!}(x-b)^k},
$$
where $f^{(k)}(b-)=\lim\limits_{x \rightarrow b-}{f^{(k)}(x)}$, for $k=0,1,\ldots,n$.

\break

Also, for $n \!\in\! \mathbb{N}$, the following functions:
$$
R_{n}^{f,\,a+}(x)
=
f(x) - T_{n-1}^{f,\,a+}(x)
$$
and
$$
R_{n}^{f,\,b-}(x)
=
f(x) - T_{n-1}^{f,\,b-}(x)
$$
are called the {\em remainder of the first {\sc Taylor}'s approximation in the right neighborhood of $a$}, and
the {\em remainder of the first {\sc Taylor}'s approximation in the left neighborhood of $b$}, respectively.

\medskip
Polynomials:
$$
\mbox{$\mathbbmsl{T}$}_n^{f;\,a+,\,b-}(x)
=
\left\{
\begin{array}{ccl}
T_{n-1}^{f,\,a+}(x)
+
\displaystyle\frac{1}{(b - a)^n}R_{n}^{f,\,a+}(b-)(x-a)^n                       &:& n \geq 1 \\[2.5 ex]
f(b-)                                                                           &:& n = 0,
\end{array}
\right.
$$
and
$$
\mbox{$\mathbbmsl{T}$}_n^{f;\,b-,\,a+}(x)
=
\left\{
\begin{array}{ccl}
T_{n-1}^{f,\,b-}(x)
+
\displaystyle\frac{1}{(a - b)^n}R_{n}^{f,\,b-}(a+)(x-b)^n                     &:& n \geq 1 \\[2.5 ex]
f(a+)                                                                           &:& n = 0,
\end{array}
\right.
$$
are called the {\em second {\sc Taylor}'s approximation in the right neighborhood of $a$}, and the {\em second {\sc Taylor}'s approximation in the left neighborhood of $b$}, respectively,  $n \!\in\! \mathbb{N}_{0}$.

\medskip

Theorem 2 in \cite{S_Wu_L_Debnath_2009} provides an important result regarding {\sc Taylor}'s approximations. We cite it below:



\begin{theorem}
\label{Theorem_1}
{\em Suppose that $f(x)$ is a real function on $(a,b)$, and that $n$ is a positive integer such that $f ^{(k)}(a+), f^{(k)}(b-)$,
for $k \!\in\! \{0,1,2, \ldots ,n\}$, exist.

\medskip
\noindent
{\boldmath $(i)$} Supposing that $(-1)^{(n)} f^{(n)}(x)$~is~in\-cre\-asing on $(a,b)$, then
for all $x \in (a,b)$ the following inequality holds$\,:$
\begin{equation}
\label{(1)}
\mbox{$\mathbbmsl{T}$}_n^{f;\,b-,\,a+}(x) < f(x) < T_n^{f,\,b-}(x).
\end{equation}
Furthermore, if $(-1)^{n} f^{(n)}(x)$ is decreasing on $(a,b)$, then the reversed inequality of {\rm (\ref{(1)})} holds.

\medskip
\noindent
{\boldmath $(ii)$} Supposing that $f^{(n)}(x)$ is increasing on $(a,b)$, then for all $x \!\in\! (a,b)$
the following inequality also holds$\,:$
\begin{equation}
\label{(2)}
\mbox{$\mathbbmsl{T}$}_n^{f;\,a+,\,b-}(x) > f(x) > T_n^{f,\,a+}(x).
\end{equation}
Furthermore, if $f^{(n)}(x)$ is decreasing on $(a,b)$, then the reversed inequality~of~\mbox{\rm (\ref{(2)})} holds.}
\end{theorem}



Let us name this theorem the {\em Theorem on double-sided \textsc{Taylor}'s approximations}. In papers
\cite{B_Malesevic_T_Lutovac_M_Rasajski_C_Mortici_Adv._Difference_Equ._2018},
\cite{M_Nenenzic_L_Zhu_AADM_2018},
\cite{M_Rasajski_T_Lutovac_B_Malesevic_JNSA_2018},
\cite{M_Rasajski_T_Lutovac_B_Malesevic_JIA_2018} and
\cite{T_Lutovac_B_Malesevic_M_Rasajaki_Results_2018}
Theorem \ref{Theorem_1} was denoted by Theorem WD.
Let us note that the proof of Theorem \ref{Theorem_1} (Theorem 2 in \cite{S_Wu_L_Debnath_2009})
was based on the  {\sc L'Hospital}'s rule for monotonicity. The same method was used in proofs of
some theorems in \cite{S_Wu_HM_Srivastva_2008a}, \cite{S_Wu_L_Debnath_2008} and \cite{S_Wu_HM_Srivastva_2008b},
that had been published earlier.

\medskip
Here, we cite a theorem (Theorem 1.1. from \cite{M_Rasajski_T_Lutovac_B_Malesevic_JNSA_2018}) that represents
a natural extension of Theorem \ref{Theorem_1} over the set of real analytic functions.

\begin{theorem}
\label{Natural_Extension_Theorem}
For the function $f:(a,b) \longrightarrow \mathbb{R}$ let there exist the power series expansion$:$
\begin{equation}
f(x)
=
\displaystyle\sum_{k=0}^{\infty}{c_{k}(x-a)^k},
\end{equation}
for every $x \!\in\! (a,b)$, where $\{c_{k}\}_{k \in \mathbb{N}_0}$ is the sequence of coefficients such that there is only
a finite number of negative coefficients, and their indices are all in the set $J \!=\! \{j_0,\ldots,j_\ell\}$.

\noindent
Then, for the function
\begin{equation}
F(x)
=
f(x)-\displaystyle\sum_{i=0}^{\ell}{c_{j_i}(x-a)^{j_i}}
=
\displaystyle\sum_{k \in N_0 \backslash\!\;\! J}{c_{k}(x-a)^k},
\end{equation}
and the sequence $\{C_{k}\}_{k \in \mathbb{N}_0}$ of the non-negative coefficients defined by
\begin{equation}
C_{k}
=
\left\{
\begin{array}{ccc}
c_{k} &:& c_{k} > 0,               \\[1.00 ex]
0     &:& c_{k} \leq 0;
\end{array}
\right.
\end{equation}
holds that$\,:$
\begin{equation}
F(x)
=
\displaystyle\sum_{k=0}^{\infty}{C_{k}(x-a)^k},
\end{equation}
for every $x \!\in\! (a,b)$.

\smallskip
Also, $F^{(k)}(a+)= k!\,C_{k}$, for $k \!\in\! \{0,1,$ $2, \ldots ,n\}$, and the following inequalities hold$\,:$
\begin{equation}
\!
\begin{array}{c}
\displaystyle\sum_{k=0}^{n}{C_k(x-a)^{k}}  <
F(x)
<                                                                                      \\[1.5 ex]
<
\displaystyle\sum_{k=0}^{n-1}{C_k(x-a)^k}
+
\frac{1}{(b-a)^n}
{\bigg (}
F(b-)
-
\displaystyle\sum_{k=0}^{n-1}{(b-a)^kC_k}
{\bigg )} (x-a)^{n},
\end{array}
\end{equation}
i.e.$\,:$
\begin{equation}
\!\!
\begin{array}{c}
\displaystyle\sum_{k=0}^{n}{\!C_k}{(x\!-\!a)^k}
+
\displaystyle\sum_{i=0}^{\ell}{\!c_{j_i}}{(x - a)^{j_i}}< f(x) <                                           \\[1.5 ex]
<
\displaystyle\sum_{k = 0}^{n - 1}{\!C_k}{(x \!-\! a)^{k}}
\!+\!
\displaystyle\sum_{i = 0}^{\ell}{\!c_{j_i}}{(x \!-\! a)^{j_i}}
\!+\!
\displaystyle\frac{{(x \!-\! a)}^n}{{(b \!-\! a)}^n}
\!\left(\!{f(b\mbox{\footnotesize $-$})
\!-\!
\displaystyle\sum_{k=0}^{n-1}{\!C_k}{{(b\!-\!a)}^k}
\!-\!
\displaystyle\sum_{i=0}^{\ell}{\!c_{j_i}}{(b\!-\!a)}^{j_i}}\!\right) ,
\end{array}
\end{equation}
for every $x \!\in\! (a,b)$.
\end{theorem}

\begin{corollary}
Let there hold the conditions from the previous theorem. If
\begin{equation}
n > \max\{j_0,\ldots,j_\ell\},
\end{equation}
then the following holds$\,:$
\begin{equation}
\!
\begin{array}{c}
\displaystyle\sum_{k=0}^{n}{\!c_k(x-a)^{k}}< f(x)<                                                                                      \\[1.5 ex]
<\displaystyle\sum_{k=0}^{n-1}{\!c_k(x-a)^k}
+
\frac{1}{(b-a)^n}
{\bigg (}
f(b-)
-
\displaystyle\sum_{k=0}^{n-1}{\!c_k(b-a)^k}
{\bigg )}(x-a)^{n},
\end{array}
\end{equation}
for every $x \!\in\! (a,b)$.
\end{corollary}

%
%

\section{Some new results on double-sided \textsc{Taylor}'s approximations }

Consider a real function $f : (a, b) \longrightarrow \mathbb{R}$ such that there exist its first and second {\sc Taylor}'s approximations on both sides,  for some $n \in \mathbb{N}$.
Let us recall the remainders in {\sc Lagrange} and the integral form, respectively,  \cite{L_E_Persson_H_Rafeiro_P_Wall_2017}:

$$
R_{n}^{f,\,a+}(x)
=
\displaystyle\frac{f^{(n)}(\xi_{a,x})}{n!}(x-a)^n ,
$$
for some $\xi_{a,x} \in (a,x)$, and 

$$
R_{n}^{f,\,a+}(x)
=
\displaystyle\frac{(x-a)^n}{(n-1)!}\displaystyle\int_{0}^{1}{f^{(n)}(a+(x-a)t) (1-t)^{n-1} \, dt}.
$$

\vspace*{-2.5 mm}

\subsection{A new proof of Theorem {\ref{Theorem_1}}}

We consider the case when $f^{(n)}(x)$ is a monotonically increasing function on $(a,b)$ for some $n \in \mathbb{N}$.
Other cases from Theorem  {\ref{Theorem_1}} are proved similarly.

\smallskip
From the {\sc Lagrange} form of the remainder and monotonicity of $f^{(n)}(x)$ on $(a,b)$ we get:
$$
\frac{f^{(n)}(a+)}{n!} < \frac{f^{(n)}(\xi_{a,x})}{n!} = \frac{f(x)-T_{n-1}^{f,\,a+}(x)}{(x-a)^n}
\quad \Longrightarrow \quad
T_{n}^{f,\,a+}(x) < f(x).
$$
since $\xi_{a,x} \in (a,x)$ for all $x \in (a,b)$.

\smallskip
Using the integral form of the remainder we obtain the following inequality for all $x \in (a,b)$:
$$
\begin{array}{rcl}
R_{n}^{f,\,a+}(x)
& = &
\displaystyle\frac{(x-a)^n}{(n-1)!}\int_{0}^{1}{f^{(n)}(a+(x-a)t) (1-t)^{n-1} \, dt}                          \\[2.5 ex]
& < &
\displaystyle\frac{(x-a)^n}{(n-1)!}\int_{0}^{1}{f^{(n)}(a+(b-a)t) (1-t)^{n-1} \, dt}                          \\[2.5 ex]
& = &
\displaystyle\frac{(x-a)^n}{(b-a)^n} R_{n}^{f,\,a+}(b-)
\qquad\qquad \Longrightarrow \qquad\qquad
f(x)
<
\mbox{$\mathbbmsl{T}$}_n^{f;\,a+,\,b-}(x).
\end{array}
$$
This completes the proof. \hfill $\Box$

\vspace*{-2.5 mm}

\subsection {Monotonicity of double-sided \textsc{Taylor}'s approximations}

\begin{proposition}
\label{Proposition_1}
Consider a real function $f : (a, b) \longrightarrow \mathbb{R}$ such that there exist its first and second {\sc Taylor}'s approximations on both sides, for some $n \in {N}_0$. Then,
\begin{equation}
\label{sgn}
{\mathop{\rm sgn}} {\Big (} \mbox{$\mathbbmsl{T}$}_{n}^{f,\, a+, \, b-}(x)\,-\,\mbox{$\mathbbmsl{T}$}_{n+1}^{f,\, a+, \, b-}(x) {\Big )}
=
{\mathop{\rm sgn}} {\Big (} {f(b-)\,-\,T_{n}^{f,\,a+}(b)} {\Big )},
\end{equation}
for all $ x\in(a,b)$.
\end{proposition}

\begin{proof}
From the definitions of the first and second {\sc Taylor}'s approximations we have:
$$ 
\begin{array}{rcl}
\mbox{$\mathbbmsl{T}$}_{n+1}^{f, a+, b-}\!\!(x)
&=&
T_{n}^{f,\,a+}\!\!(x)
+
\left(\mbox{\small $\displaystyle\frac{x-a}{b-a}$}\right)^{\!n+1} \!\!\!\!\cdot \left(f(b-)-T_n^{f,\,a}(b) \right)            \\[2.0 ex]
&=&
T_{n-1}^{f,\,a+}(x)
+
\mbox{\small $\displaystyle\frac{f^{(n)}(a+)}{n!}$}(x\!-\!a)^n
+
\left(\!
\mbox{\small $\displaystyle\frac{x-a}{b-a}$}\!\right)^{\!n} \!\!{\Big (}\mbox{\small $\displaystyle\frac{x-a}{b-a}$}\!-\!1\!+\!1\!{\Big )}
{\Big (}f(b-)-T_n^{f,\,a+}(b) {\Big )}                                                                                         \\[2.0 ex]
&=&
T_{n-1}^{f,\,a+}(x)
+
\left(\!\mbox{\small $\displaystyle\frac{x-a}{b-a}$}\!\right)^{\!n}
{\Big (}f(b-)-T_{n-1}^{f,\,a+}(b)- \mbox{\small $\displaystyle\frac{f^{(n)}(a+)}{n!}$}(b\!-\!a)^n {\Big )}+                    \\[1.0 ex]

& &
+
\mbox{\small $\displaystyle\frac{f^{(n)}(a+)}{n!}$}(x\!-\!a)^n
+
\left(\!\mbox{\small $\displaystyle\frac{x-a}{b-a}$}\!\right)^{\!n} \!\!{\Big (}\mbox{\small $\displaystyle\frac{x-a}{b-a}$}\!-\!1\!{\Big )}
{\Big (}f(b-)-T_n^{f,\,a+}(b) {\Big )} \\[2.0 ex]
&=&
\mbox{$\mathbbmsl{T}$}_{n}^{f,\, a+, \, b-}(x)
-
\mbox{\small $\displaystyle\frac{b-x}{b-a}$}\left(\!\mbox{\small $\displaystyle\frac{x-a}{b-a}$}\!\right)^{\!n}
{\Big (}f(b-)-T_n^{f,\,a+}(b) {\Big )}.
\end{array}
$$ 
Thus we have:
\begin{equation}
\mbox{$\mathbbmsl{T}$}_{n}^{f, a+, b-}\!\!(x)
-
\mbox{$\mathbbmsl{T}$}_{n+1}^{f,\, a+, \, b-}(x)
=
\frac{b-x}{b-a}\left(\!\frac{x-a}{b-a}\!\right)^{\!n}
{\Big (}f(b-)-T_n^{f,\,a+}(b) {\Big )},
\end{equation}
and the equality (\ref{sgn}) immediately follows. \hfill $\Box$
\end{proof}

Now, let us notice that if the real analytic function $f : (a, b) \longrightarrow \mathbb{R}$ satisfies the condition
$
(\forall n \in \mathbb{N}_{0})\,f^{(n)}(a+) \geq 0,
$
then, from Proposition \ref{Proposition_1} directly follows:
$$
(\forall n \in {N}_0)(\forall x \in (a,b))\,\mbox{$\mathbbmsl{T}$}_{n}^{f,\, a+, \, b-}(x)>\mbox{$\mathbbmsl{T}$}_{n+1}^{f,\, a+, \, b-}(x).
$$

\medskip
\begin{theorem}
\label{Theorem_2}
Consider a real function $f : (a, b) \longrightarrow \mathbb{R}$ such that the derivatives $f ^{(k)}(a+)$,
$k \!\in\! \{0,1,2, \ldots ,n+1\}$ exist, for some $n \in \mathbb{N}$.
\renewcommand{\labelenumi}{\textit{(\roman{enumi})}}

\smallskip
\noindent
Suppose that $f^{(n)}(x)$ and $f^{(n+1)}(x)$ are increasing on $(a,b)$, then for all $x \!\in\! (a,b)$
the following inequalities  hold$\,:$
\begin{equation}
T_n^{f,\,a+}(x)
<
T_{n+1}^{f,\,a+}(x)
<
f(x)
<
\mbox{$\mathbbmsl{T}$}_{n+1}^{f;\,a+,\,b-}(x)
<
\mbox{$\mathbbmsl{T}$}_n^{f;\,a+,\,b-}(x)
,
\end{equation}
for all $x \in (a,b)$. If $f^{(n)}(x)$ and $f^{(n+1)}(x)$ are decreasing on $(a,b)$, then for all $x \!\in\! (a,b)$
the reversed inequalities  hold.
\end{theorem}

\bigskip

\bigskip
\noindent
{\bf Case of the real analytic functions}

\bigskip
\noindent
In applications, of special interest are the real analytic functions.
\begin{theorem}
\label{Theorem_3}
Consider the real analytic functions  $f : (a, b) \longrightarrow \mathbb{R}$:
\begin{equation}
\label{f(x)}
f(x) = \sum_{k=0}^{\infty}{c_k(x-a)^k},
\end{equation}
where $c_k \in \mathbb{R}$ and  $c_k \geq 0$ for all $k \in \mathbb{N}_0$. Then,
 \begin{equation}
 \begin{array}{c}
 T_0^{f,\,a+}(x) \leq \ldots\leq T_n^{f,\,a+}(x) \leq T_{n+1}^{f,\,a+}(x) \leq \ldots                     \\[1.0 ex]
 \ldots \leq  f(x) \leq \ldots                                                         \\[1.0 ex]
 \ldots \leq  \mbox{$\mathbbmsl{T}$}_{n+1}^{f;\,a+,\,b-}(x)
 \leq   \mbox{$\mathbbmsl{T}$}_n^{f;\,a+,\,b-}(x)
 \leq  \ldots  \leq \mbox{$\mathbbmsl{T}$}_0^{f;\,a+,\,b-}(x),
 \end{array}
 \end{equation}
for all $x \in (a,b)$. If  $c_k \in \mathbb{R}$ and  $c_k \leq 0\,$ for all $k \in \mathbb{N}_0$, then the reversed inequalities hold.
\end{theorem}

%
%

\section{An application of the Theorem on double-sided \textsc{Taylor}'s approximations}

In this section we discuss an implementation of the Theorem on double-sided \textsc{Taylor}'s approximations applied to the sequence of functions:
\begin{equation}
h_{n}(x) = \displaystyle\frac{\tan x - T_{2n-1}^{\tan,\,0}(x)}{x^{2n} \tan x}
: \left(0,\mbox{\small $\displaystyle\frac{\pi}{2}$}\right) \longrightarrow R,
\end{equation}
for $n \in \mathbb{N}$.
This  sequence of functions was considered in papers \cite{Chen_Qi_2003}, \cite{Zhao_Luo_Guo_Qi_2012}.
In order to obtain estimates of functions $h_{n}(x)$, we use the well-known series expansions:

\begin{equation}
\label{tan_{series}}
\tan x
=
\displaystyle\sum_{i=1}^{\infty}{
\mbox{\small $\displaystyle\frac{ 2^{2i} (2^{2i}-1) | \textit{\textbf{B}}_{2i} | }{ (2i)! }$}x^{2i-1}},
\end{equation}
where $|x| < \mbox{\small $\displaystyle\frac{\pi}{2}$}$ and $\textit{\textbf{B}}_{k}$ is the $k$-th {\sc Bernoulli} number. Then:
\begin{equation}
T_{2n-1}^{\tan,\,0}(x)
=
\displaystyle\sum_{i=1}^{n}{
\mbox{\small $\displaystyle\frac{ 2^{2i} (2^{2i}-1) | \textit{\textbf{B}}_{2i} | }{ (2i)! }$}x^{2i-1}},
\end{equation}
for $x \in \left(0,\mbox{\small $\displaystyle\frac{\pi}{2}$}\right)$.
The main results on the functions $h_{n}(x)$,presented in the paper \cite{Zhao_Luo_Guo_Qi_2012} (see also \cite{Chen_Qi_2003}), are cited below in the following two theorems.
\begin{theorem}
\label{Theorem_4}
For $x \in \left(0,\mbox{\small $\displaystyle\frac{\pi}{2}$}\right)$ and $n \in \mathbb{N}$, we have$:$
\begin{equation}
h_{n}(x)
=
\displaystyle\sum_{j=1}^{n}{
\mbox{\small $\displaystyle\frac{2^{2(n-j+1)}(2^{2(n-j+1)}-1)|\textit{\textbf{B}}_{2(n-j+1)}|}{(2(n-j+1))!}$}
\displaystyle\sum_{k=j}^{\infty}{
\mbox{\small $\displaystyle\frac{2^{2k}|\textit{\textbf{B}}_{2k}|}{(2k)!}$} x^{2(k-j)}}}.
\end{equation}
\end{theorem}
\begin{theorem}
\label{Theorem_5}
For $x \in \left(0,\mbox{\small $\displaystyle\frac{\pi}{2}$}\right)$ and $n \in \mathbb{N}$, we have:
\begin{equation}
\label{13}
\mbox{\small $\displaystyle\frac{ 2^{2(n+1)} (2^{2(n+1)}-1) | \textit{\textbf{B}}_{2(n+1)} | }{ (2n+2)! }$}
<
h_{n}(x)
<
\mbox{\small $\left(\displaystyle\frac{2}{\pi}\right)^{\!\!2n}$},
\end{equation}
where the scalars \mbox{\small $\displaystyle\frac{ 2^{2(n+1)} (2^{2(n+1)}-1) | \textit{\textbf{B}}_{2(n+1)} | }{ (2n+2)! }$}
and \mbox{\small $\left(\displaystyle\frac{2}{\pi}\right)^{\!\!2n}$} in {\rm (\ref{13})} are the best possible.
\end{theorem}

From Theorem \ref{Theorem_4}, using the change of variables and some
algebraic transformations, immediately follows the next theorem.

\begin{theorem}
\label{Theorem_6}
For $x \in \left(0,\mbox{\small $\displaystyle\frac{\pi}{2}$}\right)$ and $n \in \mathbb{N}$, functions $h_{n}(x)$
are real analytic functions and have the following {\sc Taylor} series expansions:
\begin{equation}
h_{n}(x)
=
\displaystyle\sum_{i=0}^{\infty}{\displaystyle\sum_{j=1}^{n}{
\mbox{\small $\displaystyle\frac{2^{2(n+i+1)}(2^{2(n-j+1)}-1) \, | \textit{\textbf{B}}_{2(n-j+1)} | \, | \textit{\textbf{B}}_{2(i+j)} |}{
(2(n-j+1))! \, (2(i+j))!}$}\,x^{2i}}}.
\end{equation}
\end{theorem}

\medskip
Let us notice that the {\sc Taylor} series expansions of the functions $h_n(x)$ satisfy the conditions of Theorem \ref{Theorem_3}.

\medskip
Thus, we get the improvement of the results of Theorem \ref{Theorem_5}:

\begin{theorem}
\label{Theorem_7}
For $x \in \left(0,\mbox{\small $\displaystyle\frac{\pi}{2}$}\right)$ and $n \in \mathbb{N}$, we have
\begin{equation}
\begin{array}{c}
T_0^{\,h_n(x),\,0+}(x)
=
\mbox{\small $\displaystyle\frac{ 2^{2(n+1)} (2^{2(n+1)}-1) | \textit{\textbf{B}}_{2(n+1)} | }{ (2n+2)! }$} <                  \\[2.5 ex]
< T_2^{\,h_n(x),\,0+}(x)
< \ldots
< T_{2m}^{\,h_n(x),\,0+}(x)
< T_{2m+2}^{\,h_n(x),\,0+}(x)
< \ldots                                                                                                                       \\[2.5 ex]
\ldots<
h_{n}(x)
<\ldots                                                                                                                              \\[1.5 ex]
 \ldots < \mbox{$\mathbbmsl{T}$}_{2m+2}^{\,h_n(x);\,0+,\,\frac{\pi}{2}-}(x)
< \mbox{$\mathbbmsl{T}$}_{2m}^{\,h_n(x);\,0+,\,\frac{\pi}{2}-}(x)
< \ldots
< \mbox{$\mathbbmsl{T}$}_2^{\,h_n(x);\,0+,\,\frac{\pi}{2}-}(x) <                                                               \\[2.5 ex]
< \mbox{$\mathbbmsl{T}$}_0^{\,h_n(x);\,0+,\,\frac{\pi}{2}-}(x)
=
\mbox{\small $\left(\displaystyle\frac{2}{\pi}\right)^{\!\!2n}$}.
\end{array}
\end{equation}
\end{theorem}

%
%

\section{More examples of double-sided Taylor's approximations}

In this section we give two examples of some analytic inequalities recently proved using the results of
Theorem \ref{Theorem_1}. Also, we illustrate the application of double-sided {\sc Taylor} approximations and
Theorem \ref{Theorem_3} in the generalizations and improvements of some analytic inequalities.

\bigskip
\noindent    
{\bf Example 1.}
In \cite{L_Debnath_C_Mortici_L_Zhu_2015} the following improvement of {\sc Ste\v ckin}'s inequality,
in the left neighborhood of $b=\mbox{\small $\displaystyle\frac{\pi}{2}$}$,  was proposed and proved:
\begin{equation}
\label{f_23}
Q_{1}(x)
=
\mbox{\small $\displaystyle\frac{2}{\pi}$}
-
\mbox{\small $\displaystyle\frac{1}{2}$}\left(\mbox{\small $\displaystyle\frac{\pi}{2}$} - x\right)
<
\tan \, x - \mbox{\small $\displaystyle\frac{4x}{\pi(2\pi-x)}$}
<
\mbox{\small $\displaystyle\frac{2}{\pi}$}
-
\mbox{\small $\displaystyle\frac{1}{3}$}\left(\mbox{\small $\displaystyle\frac{\pi}{2}$} - x\right)
=
R_{1}(x),
\end{equation}
for  $x \!\in\! \left(0, \mbox{\small $\displaystyle\frac{\pi}{2}$}\right)$.
In  \cite{M_Nenenzic_L_Zhu_AADM_2018} the  inequality~(\ref{f_23}) was further generalized.
The starting point   was the following real function:
\begin{equation}
g(t)
\!=\!
\cot \, t - \frac{1}{t} + \frac{2}{\pi}
:
\left(0, \mbox{\small $\displaystyle\frac{\pi}{2}$} \right) \longrightarrow R,
\end{equation}
for which it is fulfilled
\begin{equation}
\label{g-and-f}
g\left(\mbox{\small $\displaystyle\frac{\pi}{2}$}-x\right) = \tan \, x - \mbox{\small $\displaystyle\frac{4x}{\pi(2\pi-x)}$},
\end{equation}
for  $x \!\in\! \left(0, \mbox{\small $\displaystyle\frac{\pi}{2}$}\right)$.
It has been shown that the function $g(t)$  satisfies the conditions of Theorem~\ref{Theorem_1}.
Namely, it has the following power series expansion
\begin{equation}
g(t)
\!=\!
\mbox{\small $\displaystyle\frac{2}{\pi}$}
-
\displaystyle\sum \limits_{k = 1}^{\infty} \frac{2^{2k} {\left| {{\mbox{\textit{\textbf{B}}}_{2k}}} \right| }}{{ (2k)! }}{t^{2k-1}}
\end{equation}
which converges for $t \in \left(0,\mbox{\small $\displaystyle\frac{\pi}{2}$}\right)$,
and it is true
$$
\mbox{$g(0+)=\!\lim\limits_{t \rightarrow 0+}{\!g(t)}
=
\mbox{\small $\displaystyle\frac{2}{\pi}$}$}
\quad\mbox{and}\quad
\mbox{$g\left(\mbox{\small $\displaystyle\frac{\pi}{2}$}-\right)=\!
\lim\limits_{t \rightarrow \pi/2-}{\!g(t)}=0$}.
$$
The function $g(t)$ also satisfies the conditions of Theorem~\ref{Theorem_3}. Based on this, the following result was proposed in
\cite{M_Nenenzic_L_Zhu_AADM_2018} (Theorem 3) for the function $f(x)=g\left(\mbox{\small $\displaystyle\frac{\pi}{2}$}-x\right)$:
\begin{theorem}
For every $x \!\in\! \left(0, \mbox{\small $\displaystyle\frac{\pi}{2}$}\right)$
and $m \in N$, $m \ge 2$, the following inequalities hold$:$

\begin{equation}
\mbox{$\mathbbmsl{T}$}_{2m-1}^{g;\,0+,\,\pi/2-}\!\left(\mbox{\small $\displaystyle\frac{\pi}{2}$}-x\right)\, <  f(x) \, <   T_{2m-1}^{g,\;0}\left(\mbox{\small $\displaystyle\frac{\pi}{2}$}-x\right),
\end{equation}
where
\begin{equation}
\begin{array}{c}
\mbox{$\mathbbmsl{T}$}_{2m-1}^{g;\,0+,\,\pi/2-}\!\left(\mbox{\small $\displaystyle\frac{\pi}{2}$}-x\right)=       \\[2.0ex]

 =
\mbox{\small $\displaystyle\frac{2}{\pi}$}
-
\displaystyle\sum \limits_{k = 1}^{m-1} \frac{2^{2k} {\left| {{\mbox{\textit{\textbf{B}}}_{2k}}} \right| }}{{(2k)!}}{
\left(\mbox{\small $\displaystyle\frac{\pi}{2}$}-x\right)^{\!2k-1}}
\!\!+
{\displaystyle\sum \limits_{k = 1}^{m-1} \frac{2^{2k} {\left| {{\mbox{\textit{\textbf{B}}}_{2k}}} \right| }}{{ (2k)! }}
\left( \mbox{\small $\displaystyle\frac{2}{\pi}$} \right)^{2(m+k-1)}}
\!\left(\mbox{\small $\displaystyle\frac{\pi}{2}$}-x\right)^{\!2m-1}
\end{array}
\end{equation}
and
\begin{equation}
T_{2m-1}^{g,\;0}\left(\mbox{\small $\displaystyle\frac{\pi}{2}$}-x\right)
=
\mbox{\small $\displaystyle\frac{2}{\pi}$}
-
\displaystyle\sum \limits_{k = 1}^{m} \frac{2^{2k} {\left| {{\mbox{\textit{\textbf{B}}}_{2k}}} \right| }}{{ (2k)! }}{
\left(\mbox{\small $\displaystyle\frac{2}{\pi}$}-x\right)^{\!2k-1}}.
\end{equation}

\end{theorem}

It is easy to check that the function $g(t)$ also satisfies the conditions of Theorem~\ref{Theorem_3}.
Therefore, for the function $f(x)=g\left(\mbox{\small $\displaystyle\frac{\pi}{2}$}-x\right)$
the following assertion directly follows:

\bigskip

\begin{theorem}
\label{unlimited}
For every $x \!\in\! \left(0, \mbox{\small $\displaystyle\frac{\pi}{2}$}\right)$ and $m \in N$, $m \ge 2$, the  following inequalities hold$:$
 \begin{equation}
 \begin{array}{c}
 T_1^{g,\,0+}\left(\mbox{\small $\displaystyle\frac{\pi}{2}$}-x\right) \leq \ldots\leq T_{2m-1}^{g,\,0+}\left(\mbox{\small $\displaystyle\frac{\pi}{2}$}-x\right) \leq T_{2m+1}^{g,\,0+}\left(\mbox{\small $\displaystyle\frac{\pi}{2}$}-x\right) \leq \ldots                     \\[2.0 ex]
 \ldots \leq  f(x) \leq \ldots                                                         \\[2.0 ex]
 \ldots \leq  \mbox{$\mathbbmsl{T}$}_{2m+1}^{g;\,0+,\,\frac{\pi}{2}-}\left(\mbox{\small $\displaystyle\frac{\pi}{2}$}-x\right)
 \leq   \mbox{$\mathbbmsl{T}$}_{2m-1}^{g;\,0+,\,\frac{\pi}{2}-}\left(\mbox{\small $\displaystyle\frac{\pi}{2}$}-x\right)
 \leq  \ldots  \leq \mbox{$\mathbbmsl{T}$}_1^{g;\,0+,\,\frac{\pi}{2}-}\left(\mbox{\small $\displaystyle\frac{\pi}{2}$}-x\right).
 \end{array}
 \end{equation}
\end{theorem}

\bigskip

Finally, from the previous two theorems an improvement of the inequality (\ref{f_23}) directly follows.
For example, if $m=1$, we have:
$$
\begin{array}{rcl}
\mbox{$\mathbbmsl{T}$}_{1}^{g;\,0+,\,\pi/2-}\left(\mbox{\small $\displaystyle\frac{\pi}{2}$}-x\right)
& = &
\mbox{\small $\displaystyle\frac{2}{\pi}$}
-
\mbox{\small $\displaystyle\frac{4}{\pi^2}$}\left(\mbox{\small $\displaystyle\frac{\pi}{2}$} - x\right) \, \leq         \\[2.0 ex]
& \leq &
\tan \, x - \mbox{\small $\displaystyle\frac{4x}{\pi(2\pi-x)}$} \, \leq                                                 \\[2.0 ex]
& \leq &
\mbox{\small $\displaystyle\frac{2}{\pi}$}
-
\mbox{\small $\displaystyle\frac{1}{3}$}\left(\mbox{\small $\displaystyle\frac{\pi}{2}$} - x\right)
=
T_{1}^{g,\,0}\left(\mbox{\small $\displaystyle\frac{\pi}{2}$} - x\right),
\end{array}
$$
which further implies the following:
$$
Q_{1}(x)
<
\mbox{$\mathbbmsl{T}$}_{1}^{g;\,0+,\,\pi/2-}\left(\mbox{\small $\displaystyle\frac{\pi}{2}$}-x\right)
\leq
\tan \, x - \mbox{\small $\displaystyle\frac{4x}{\pi(2\pi-x)}$}
\leq
T_{1}^{g,\,0+}\left(\mbox{\small $\displaystyle\frac{\pi}{2}$}-x\right)
=
R_{1}(x),
$$
for $x \!\in\! \left(0, \mbox{\small $\displaystyle\frac{\pi}{2}$}\right)$.

\medskip
Note that the same approach (based on Theorem \ref{Theorem_1} and Theorem \ref{Theorem_3}) enables generalizations
of the inequalities from \cite{M_Nenenzic_L_Zhu_AADM_2018} connected with the function
$$
f(x) = \left(\pi^2 - 4 x^2 \right) \frac{\tan x}{x}
:
\left(0, \mbox{\small $\displaystyle\frac{\pi}{2}$} \right) \longrightarrow R.
$$

\bigskip

\noindent    
{\bf Example 2.}
In {\rm \cite{C_Mortici_2011}} (Theorem 5)  the following inequality was proved:
\begin{equation}
\label{Theorem_5_C.Mortici}
\displaystyle
  2 + \frac{2}{45}\,x^4 <  \left(\frac{x}{\sin{x}}\right)^{\!2} + \frac{x}{\tan{x}} ~~~~ \mbox{\rm for} ~~~
0<x<\pi/2.
\end{equation}
In order to refine the previous inequality,  the following real function was considered in
\cite{B_Malesevic_T_Lutovac_M_Rasajski_C_Mortici_Adv._Difference_Equ._2018}:
$$
\displaystyle f(x) = \left(\frac{x}{\sin{x}}\right)^{\!2} + \frac{x}{\tan{x}} ~~~ \mbox{for} ~~~
0<x<\pi/2.
$$
It has been shown that the above function satisfies the conditions of Theorem \ref{Theorem_1}. \\

Namely, it  has the following power series expansion
\begin{equation}
f(x)
\!=\!
2 + \sum \limits_{k = 2}^{\infty} \frac{{\left| {{\mbox{\textit{\textbf{B}}}_{2k}}} \right| (2k - 2){4^k}}}{{(2k)!}}{x^{2k}},
\end{equation}
which converges for $x \in \left(0,\mbox{\small $\displaystyle\frac{\pi}{2}$}\right)$, and it is true
$$
\mbox{$f(0+)=\!\lim\limits_{x \rightarrow 0+}{\!f(x)}=2$}
\quad\mbox{and}\quad
\mbox{$f\left(\mbox{\small $\displaystyle\frac{\pi}{2}$}-\right)=\!
\lim\limits_{x \rightarrow \pi/2-}{\!f(x)}=\mbox{\small $\displaystyle\frac{\pi^2}{4}$}.$}
$$
Based on this, the following result was proposed and proved in \cite{B_Malesevic_T_Lutovac_M_Rasajski_C_Mortici_Adv._Difference_Equ._2018}
(Theorem 5):

\begin{theorem}
$\!\!$ For every $x \!\in\! \left(0, \mbox{\small $\displaystyle\frac{\pi}{2}$}\right)$ and $m \in N$, $m \ge
2$, the following inequalities hold:
\begin{equation}
\label{ineq_AIDE}
T_{2m}^{f,\;0}(x) <   f(x) < \mbox{$\mathbbmsl{T}$}_{2m}^{f;\,0+,\,\pi/2-}(x),
\end{equation}
where

\begin{equation}
T_{2m}^{f,\;0}(x)
=
2
+
\!\displaystyle\sum \limits_{k = 2}^m
\mbox{\small $\displaystyle\frac{{\left| {{\textit{\textbf{B}}_{2k}}} \right|(2k - 2){4^k}}}{{(2k)!}}$}{x^{2k}}
\end{equation}
and
\begin{equation}
\begin{array}{c}
\mbox{$\mathbbmsl{T}$}_{2m}^{f;\,0+,\,\pi/2-}(x)
=\\[2.0ex]
=
\displaystyle
2 + \displaystyle\sum \limits_{k = 2}^{m - 1}
\mbox{\small $\displaystyle\frac{{\left| {{\textit{\textbf{B}}_{2k}}} \right|(2k - 2){4^k}}}{{(2k)!}}$}{x^{2k}}  +
{\left( \mbox{\small $\displaystyle\frac{{2}}{\pi}$} \right)^{\!\!2m}}
\!\!\left( \! {
\mbox{\small $\displaystyle\frac{{{\pi^2}}}{4}$}
- 2 -\! \displaystyle\sum \limits_{k = 2}^{m - 1}
\mbox{\small $\displaystyle\frac{{\left| {{\textit{\textbf{B}}_{2k}}} \right|(2k - 2){4^k}}}{{(2k)!}}$}{{\left(
\mbox{\small $\displaystyle\frac{\pi}{2}$} \right)}^{\!2k}}} \! \right)\!x^{2m}.
\end{array}
\end{equation}
\end{theorem}
In \cite{B_Malesevic_T_Lutovac_M_Rasajski_C_Mortici_Adv._Difference_Equ._2018}
the polynomials \mbox{$T_{m}^{f,\,0+}\!(x)$} and $\mbox{$\mathbbmsl{T}$}_{m}^{f;\,0+,\,\pi/2-}\!(x)$ are calculated and the concrete inequalities
$$T_{m}^{f,\,0+}\!(x) < f(x) < \mbox{$\mathbbmsl{T}$}_{m}^{f;\,0+,\,\pi/2-}\!(x)$$
are given for \mbox{$x \!\in\! \left(0, \mbox{\small $\displaystyle\frac{\pi}{2}$}\right)$} and for $m = 2, 3, 4, 5$.

\medskip

It is easy to check that the function $f(x)$ also satisfies the conditions of Theorem \ref{Theorem_3},
and hence the following generalizations  of the inequality (\ref{ineq_AIDE})$\;$i.e. of the inequality (\ref{Theorem_5_C.Mortici}) are true:

\begin{theorem}
For every $x \!\in\! \left(0, \mbox{\small $\displaystyle\frac{\pi}{2}$}\right)$ and $m \in N$, $m \ge
2$, the following inequalities hold$:$
 \begin{equation}
 \begin{array}{c}
 T_0^{f,\,0+}(x) \leq \ldots\leq T_{2m}^{f,\,0+}(x) \leq T_{2m+2}^{f,\,0+}(x) \leq \ldots                     \\[1.0 ex]
 \ldots \leq  f(x) \leq \ldots                                                         \\[1.0 ex]
 \ldots \leq  \mbox{$\mathbbmsl{T}$}_{2m+2}^{f;\,0+,\,\pi/2-}(x)
 \leq   \mbox{$\mathbbmsl{T}$}_{2m}^{f;\,0+,\,\pi/2-}(x)
 \leq  \ldots  \leq \mbox{$\mathbbmsl{T}$}_0^{f;\,0+,\,\pi/2-}(x)
 \end{array}
 \end{equation}
\end{theorem}

\medskip
The same approach, based on Theorem \ref{Theorem_1} and Theorem \ref{Theorem_3}, provides generalizations of the inequalities
from \cite{B_Malesevic_T_Lutovac_M_Rasajski_C_Mortici_Adv._Difference_Equ._2018} related to the function
$$
f(x) = 3 \frac{x}{\sin x} + \cos x
:
\left(0, \mbox{\small $\displaystyle\frac{\pi}{2}$} \right) \longrightarrow R.
$$

%
%

\section{Conclusion}
Even though \textsc{Taylor}'s approximations represent a few centuries old topic, they are still present
in research nowadays in many areas of science and engineering. Let us note that many results regarding
\textsc{Taylor}'s approximations are presented in well-known monographs \cite{D_S_Mitrinovic_1970} and
\cite{Milovanovic_Rassias_2014}. Historically speaking the second {\sc Taylor}'s approximation was mentioned
in 1851 in the proof of the {\sc Taylor}'s formula with the {\sc Lagrange} remainder in the paper  \cite{H_Cox_1851}
by {\sc H. Cox}, see also \cite{L_E_Persson_H_Rafeiro_P_Wall_2017}.

\smallskip
Let us mention that in papers
\cite{Zhu_Hua_2010},
\cite{Milica_Makragic_JMI_2017},
\cite{H_Alzer_M_K_Kwong_2017},
 \cite{M_Rasajski_T_Lutovac_B_Malesevic_JNSA_2018},
\cite{M_Rasajski_T_Lutovac_B_Malesevic_JIA_2018} and
\cite{T_Lutovac_B_Malesevic_M_Rasajaki_Results_2018}
double-sided {\sc Ta\-ylo\-r}'s approximations are used to obtain corresponding inequalities.
Results of these papers can be further organized and made more precise using Theorem \ref{Theorem_3}
so we get the order among the functions occurring within these inequalities.
Similar to double-sided  {\sc Ta\-ylo\-r}'s approximations, in papers
\cite{B_Banjac_M_Nenenzic_B_Malesevic_Telfor_2015},
\cite{B_Malesevic_M_Makragic_JMI_2016},
\cite{M_Nenezic_B_Malesevic_C_Mortici_AMC_2016},
\cite{B_Banjac_M_Makragic_B_Malesevic_Results_2016},
\cite{T_Lutovac_B_Malesevic_C_Mortici_JIA_2017},
\cite{B_Malesevic_M_Rasajski_T_Lutovac_JIA_2017},
\cite{B_Malesevic_I_Jovovic_B_Banjac_JMI_2017},
\cite{B_Malesevic_T_Lutovac_B_Banjac_JMI_2018},
\cite{B_Malesevic_M_Rasajski_T_Lutovac_MPE_2018} i
\cite{B_Malesevic_T_Lutovac_B_Banjac_Acta_2017}
the finite expansions are used in the proofs of some mixed-trigonometric polynomial inequalities,
as well as in some inequalities which can be reduced to mixed-trigonometric polynomial inequalities.

\smallskip
Currently, we are working on developing a computer system for automatic proving of some classes
of analytic inequalities based on the results in the mentioned papers.

\bigskip
\textbf{Acknowledgment.}
Research of the first and second and third author was supported in part by the Serbian Ministry of
Education, Science and Technological Development, under Projects ON 174032 \& III 44006, ON 174033
and TR 32023, respectively.


\begin{thebibliography}{99}%

\bibitem{H_Cox_1851}
H.~Cox, \textit{A demonstration of Taylor's theorem}, Cambridge and Dublin Math. J., 6, (1851) 80--81.

\bibitem{D_S_Mitrinovic_1970}
D.~S.~Mitrinovi\' c, \textit{Analytic inequalities}, Springer (1970).

\bibitem{Chen_Qi_2003}
Ch.-P. Chen, F. Qi, {\em A double inequality for remainder of power series of tangent function},
Tamkang J. Math. {\bf 34}:3 (2003), 351--355.

\bibitem{S_Wu_HM_Srivastva_2008a}
S.-H. Wu, H.M. Srivastva, \textit{A further refinement of a Jordan type inequality and its applications},
Appl. Math. Comput. {\bf 197} (2008), 914--923.

\bibitem{S_Wu_L_Debnath_2008}
S.-H. Wu, L. Debnath, \textit{Jordan-type inequalities for differentiable functions and their applications},
Appl. Math. Lett. {\bf 21}:8 (2008), 803--809.

\bibitem{S_Wu_HM_Srivastva_2008b}
S.-H. Wu, H. M. Srivastava, \textit{A further refinement of Wilker's inequality},
Integral Transforms Spec. Funct. {\bf 19}:9-10 (2008), 757--765.

\bibitem{S_Wu_L_Debnath_2009}
S.~Wu, L.~Debnath, \textit{A generalization of L'Hospital-type rules for monotonicity and its application},
Appl. Math. Lett. {\bf 22}:2 (2009), 284--290.

\bibitem{Zhu_Hua_2010}
L. Zhu, J. Hua, \textit{Sharpening the Becker-Stark inequalities},
J. Inequal. Appl. {\bf 2010} (2010), 1--4.

\bibitem{C_Mortici_2011}
C. Mortici, \textit{The natural approach of Wilker-Cusa-Huygens inequalities}, Math. Inequal. Appl. {\bf 14}:3 (2011), 535--541.

\bibitem{Zhao_Luo_Guo_Qi_2012}
J.-L. Zhao, Q.-M. Luo, B.-N. Guo, and F. Qi, \textit{Remarks on inequalities for the tangent function},
Hacettepe J. Math. Statist. {\bf 41}:4 (2012), 499--506.

\bibitem{Milovanovic_Rassias_2014}
G.$\,$Milovanovi\'c,$\,$M.$\,$Rassias (ed.),
\textit{Analytic Number Theory, Approximation Theory and Special Functions}, Springer 2014.
{\big (}Chapter:~G.$\,$D. Anderson, M. Vuorinen, X. Zhang: \textit{Topics in Special Functions III},~297--345.{\big )}

\bibitem{B_Banjac_M_Nenenzic_B_Malesevic_Telfor_2015}
B.~Banjac, M.~Nenezi\'c, B.~Male\v sevi\'c, \textit{Some applications of Lambda-method for obtaining approximations in filter design},
Proceedings of 23-rd TELFOR conference, pp. 404-406, Beograd 2015.

\bibitem{L_Debnath_C_Mortici_L_Zhu_2015}
L. Debnath, C. Mortici, L. Zhu, \textit{Refinements of Jordan-Steckin and Becker-Stark inequalities},
Results Math. 67 (1-2) (2015), 207--215.

\bibitem{B_Malesevic_M_Makragic_JMI_2016}
B.~Male\v sevi\' c, M.~Makragi\' c, \textit{A Method for Proving Some Inequalities on Mixed Trigonometric Polynomial Functions},
J. Math. Inequal. {\bf 10}:3 (2016), 849--876.

\bibitem{M_Nenezic_B_Malesevic_C_Mortici_AMC_2016}
M.~Nenezi\' c, B.~Male\v sevi\' c, C.~Mortici, \textit{New approximations of some expressions involving trigonometric functions},
Appl. Math. Comput. {\bf 283} (2016), 299--315.

\bibitem{B_Banjac_M_Makragic_B_Malesevic_Results_2016}
B.~Banjac, M.~Makragi\' c, B.~Male\v sevi\' c \textit{Some notes on a method for proving inequalities by computer},
Results Math. {\bf 69}:1 (2016), 161--176.

\bibitem{L_E_Persson_H_Rafeiro_P_Wall_2017}
L.~E.~Persson, H.~Rafeiro, P.~Wall, \textit{Historical synopsis of the Taylor remainder}, Note Mat. {\bf 37}:1 (2017), 1--21.

\bibitem{Milica_Makragic_JMI_2017}
M.~Makragi\' c, \textit{A method for proving some inequalities on mixed hyper\-bolic-trigonometric polynomial functions},
J. Math. Inequal. {\bf 11}:3 (2017), 817--829.

\bibitem{T_Lutovac_B_Malesevic_C_Mortici_JIA_2017}
T.~Lutovac, B.~Male\v sevi\' c, C.~Mortici, \textit{The natural algorithmic approach of mixed trigonometric-polynomial problems},
J. Inequal. Appl. {\bf 2017}:116 (2017), 1--16.

\bibitem{B_Malesevic_M_Rasajski_T_Lutovac_JIA_2017}
B.~Male\v sevi\' c, M.~Ra\v sajski, T.~Lutovac, \textit{Refinements and generalizations of some inequalities of Shafer-Fink's type
for the inverse sine function}, J. Inequal. Appl. {\bf 2017}:275 (2017), 1--9.

\bibitem{B_Malesevic_I_Jovovic_B_Banjac_JMI_2017}
B.~Male\v sevi\' c, I.~Jovovi\' c, B.~Banjac, \textit{A proof of two conjectures of Chao-Ping Chen for inverse trigonometric functions},
J. Math. Inequal. {\bf 11} (1) (2017),  151--162.

\bibitem{H_Alzer_M_K_Kwong_2017}
H.~Alzer, M.~K.~Kwong, \textit{On Jordan's inequality}, Period. Math. Hung. {\bf 77}:2 (2018), 191--200.

\bibitem{B_Malesevic_T_Lutovac_M_Rasajski_C_Mortici_Adv._Difference_Equ._2018}
B.~Male\v sevi\' c, T.~Lutovac, M.~Ra\v sajski, C.~Mortici, \textit{Extensions of the natural approach to refinements and generalizations
of some trigonometric inequalities}, Adv. Difference Equ. {\bf 2018}:90 (2018), 1--15.

\bibitem{M_Nenenzic_L_Zhu_AADM_2018}
M.~Nenezi\' c, L.~Zhu, \textit{Some improvements of Jordan-Steckin and Becker-Stark inequalities},
Appl. Anal. Discrete Math. {\bf 12} (2018), 244--256.

\bibitem{M_Rasajski_T_Lutovac_B_Malesevic_JNSA_2018}
M.~Ra\v sajski, T.~Lutovac,~B. Male\v sevi\' c, \textit{Sharpening and generalizations of Shafer-Fink and Wilker type inequalities$:$ a new approach},
J. Nonlinear Sci. Appl. {\bf 11}:7 (2018), 885--893.

\bibitem{M_Rasajski_T_Lutovac_B_Malesevic_JIA_2018}
M.~Ra\v sajski, T.~Lutovac, B.~Male\v sevi\' c, \textit{About some exponential inequalities related to the sinc function},
J. Inequal. Appl. {\bf 2018}:150 (2018), 1--10.

\bibitem{T_Lutovac_B_Malesevic_M_Rasajaki_Results_2018}
T.~Lutovac, B.~Male\v sevi\' c, M.~Ra\v sajski, \textit{A new method for proving some inequalities related to several special functions},
Results Math. {\bf 73}:100 (2018), 1--15.

\bibitem{B_Malesevic_T_Lutovac_B_Banjac_JMI_2018}
B. Male\v sevi\' c, T. Lutovac, B. Banjac, \textit{A proof of an open problem of Yusuke Nishizawa for a power-exponential function},
J. Math. Inequal. {\bf 12}:2 (2018), 473-485.

\bibitem{B_Malesevic_M_Rasajski_T_Lutovac_MPE_2018}
B.~Male\v sevi\' c, M.~Ra\v sajski, T.~Lutovac, \textit{Refined estimates and generalizations of inequalities related
to the arctangent function and Shafer's inequality}, Math. Probl. Eng. {\bf 2018} Article ID 4178629, 1--8.

\bibitem{B_Malesevic_T_Lutovac_B_Banjac_Acta_2017}
B.~Male\v sevi\' c, T.~Lutovac, B.~Banjac, \textit{One method for proving some classes of exponential analytical inequalities},
Accepted in Filomat ACTA 2017 - Special issue (2018)


\end{thebibliography}
\end{document}